\documentclass[12pt]{imsart}

\usepackage[margin=0.9in]{geometry}
\usepackage{graphicx}
\usepackage[font=large]{caption} 

\usepackage[toc]{appendix}
\RequirePackage[OT1]{fontenc}
\RequirePackage[authoryear]{natbib}
\usepackage{adjustbox}
\usepackage{amsmath,amssymb,amsthm,bm,latexsym}
\usepackage{graphicx,psfrag,epsf}
\usepackage{epigraph,xcolor}
\usepackage{enumerate}
\usepackage{url}
\usepackage{booktabs}
\usepackage{multirow}
\usepackage{hyperref}
\hypersetup{colorlinks=true,citecolor=blue,pdfpagemode=FullScreen, linktocpage=true}
\usepackage{cleveref}
\usepackage{amsfonts}
\usepackage{amsmath}
\usepackage{orcidlink}

\startlocaldefs
\def\E{\mathbb{E}}

\def\Z{\mathbb{Z}}

\def\R{\mathbb{R}}

\newcommand{\RNum}[1]{\uppercase\expandafter{\romannumeral #1\relax}}

\DeclareMathOperator*{\argmax}{arg\,max}

\newtheorem{theorem}{Theorem}[section]

\newtheorem{definition}{Definition}[section]

\newtheorem{lemma}{Lemma}[section]
\newtheorem{remark}{Remark}[section]
\newtheorem{proposition}{Proposition}[section]

\numberwithin{equation}{section}  
\endlocaldefs

\begin{document}

\begin{frontmatter}
\title{Detecting Change-points in Mean of Multivariate Time Series}
\runtitle{Detecting Change-points in Mean of Multivariate Time Series}

\begin{aug}
\author[MCM]{\fnms{Ramkrishna Jyoti} \snm{Samanta}\,\orcidlink{0009-0003-3021-6195}\ead[label=e1]{akashnilsamanta@gmail.com}}
  
\address[MCM]{McMaster University}  
\printead{e1}

\runauthor{RJ Samanta}
\end{aug}

\begin{abstract}
This work delves into presenting a probabilistic method for analyzing linear process data with weakly dependent innovations, focusing on detecting change-points in the mean and estimating its spectral density.  We develop a test for identifying change-points in the mean of data coming from such a model, aiming to detect shifts in the underlying distribution. Additionally, we propose a consistent estimator for the spectral density of the data, contingent upon fundamental assumptions, notably the long-run variance. By leveraging probabilistic techniques, our approach provides reliable tools for understanding temporal changes in linear process data. Through theoretical analysis and empirical evaluation, we demonstrate the efficacy and consistency of our proposed methods, offering valuable insights for practitioners in various fields dealing with time series data analysis. Finally, we implemented our method on bitcoin data for identifying the time points of significant changes in its stock price.
\end{abstract}

\tableofcontents
\end{frontmatter}

\section{Introduction} 
\label{sec:introduction}
A structural break refers to a sudden and significant change in the underlying data-generating process of a time series. This change can manifest in various parameters of the model, such as the mean, variance, or autoregressive coefficients and a change-point refers to a specific time point within a time series where a structural break occurs. For example, consider a temperature dataset spanning several decades. Suppose the data exhibits a stable pattern with gradual temperature increases over time, representing a consistent warming trend. However, in recent years, there has been a sudden and persistent drop in temperatures, leading to colder conditions. The time point where this abrupt change occurs marks the change-point, indicating a structural break in the temperature trend. Detecting structural breaks has captivated researchers' interest since  \cite{pagetest} proposed a test for it as early as in 1955. Identifying changes in the mean of a multivariate data is a fundamental problem of interest in various disciplines such as signal processing, finance, environmental monitoring, and healthcare. Change-point detection methods are pivotal in pinpointing abrupt shifts or structural breaks in the underlying distribution of multivariate data over time or space. These methods are indispensable for comprehending and analyzing dynamic systems, detecting anomalies, and making well-informed decisions based on evolving data patterns. Various change-point detection methods have been discussed briefly in   \cite{johnson1998encyclopedia}, where an overview is presented of various techniques used by researchers in developing methods for change-point detection.

Multivariate change-point detection presents unique challenges compared to the univariate scenario due to the interdependence among the observed variables. In the case of multivariate data, changes in one variable may correlate with changes in others, necessitating sophisticated techniques to accurately detect shifts in the mean across multiple dimensions. A time series $\{X_t\}$ is called a linear process if it can be expressed as:

\[
X_t = \mu + \sum_{j=-\infty}^{\infty} \psi_j \varepsilon_{t-j},
\],
where $\mu$ is the mean of $X_t$, $\{\psi_j\}$ are the coefficients with $\sum_{j=-\infty}^{\infty} |\psi_j| < \infty$, and $\{\varepsilon_t\}$ is white noise with zero mean and constant variance. In certain types of data, it will be feasible to model data with some level of dependence among its coordinates. This allows for the introduction of dependence in the innovations. In our case, we specifically consider linear process data with \(m\)-dependent $\{\varepsilon_t\}$'s. A sequence of random variables $\{X_t\}_{t \in \mathbb{Z}}$ is said to be \(m\)-dependent if for any two sets of indices $I$ and $J$ such that the minimum distance between elements of $I$ and elements of $J$ is greater than $m$, the sets of random variables $\{X_i : i \in I\}$ and $\{X_j : j \in J\}$ are independent. Formally,

\[
\{X_i : i \in I\} \perp\!\!\!\perp \{X_j : j \in J\} \quad \text{whenever} \quad \min_{i \in I, j \in J} |i - j| > m
\]

Here, $\perp\!\!\!\perp$ denotes the independence of the two sets of random variables. To conduct a change-point analysis, we first need to perform hypothesis testing to determine whether a change-point exists at all. This typically involves a statistic dependent on $N$, the length of the data, whose distribution for finite $N$ is unknown. Hence, we consider data with a large length ($N\rightarrow \infty$), yielding a suitable statistic with known critical values from \cite{Kiefer}. Moreover, \cite{zeileis2002detecting} stated that if a time series is observed for a long enough period of time, there would be some economic or political or climatic factors that would cause the structure of the series to change at some points. This also entails estimating the long-run covariance matrix of the data, which leads to the estimation of the spectral density of the data.

Let $\{X_t\}$ be a stationary time series with mean $\mu$. The spectral density function $f(\lambda)$ at frequency $\lambda$ is defined as the Fourier transform of the autocovariance function $\Gamma(h)$. For a discrete-time stationary process, the spectral density is given by:

$$f(\lambda) = \frac{1}{2\pi}\sum_{h=-\infty}^{\infty} \Gamma(h)e^{-ih\lambda}, \quad -\pi \leq \lambda \leq \pi.$$
Estimating spectral density has been studied extensively for univariate linear process data with independent innovations, as can be seen in \cite{brockwell_davis}. Here, we aim to derive a consistent estimator of the spectral density for such data, especially to estimate the long-run covariance matrix of the data. This, in turn, facilitates us to provide an estimator for the change-point under the assumption that there is a change-point. As usual, we consider a $d$ dimensional time series $X_t = (X^1_t,\ldots,X^d_t)$, and $t\in\{1,\ldots,N\}$ corresponding to the time-points of observations. In our case, we have:
$$
\begin{aligned}
&H_0: \quad \E[X_1]=\ldots=\E[X_N]\\
&\text{vs}\\
&H_a: \quad \exists \hspace{0.1in} T^* \text{ such that } \E[X_1]=\ldots=\E[X_{T^*}]\neq \E[X_{T^*+1}], \text{where }T^*<N.
\end{aligned}
$$
This is the scenario when we have a single change-point in our data. There have been many works in the literature discussing various types of change-point detection methods in multivariate data. Some of them include \cite{matteson2014nonparametric}, \cite{kuncheva2011change}, \cite{gao2020change}, \cite{horvath2012change}, etc. 

Section \ref{pre} of this work provides an overview of the necessary preliminaries. In Section \ref{theory}, we develop the theoretical framework required for hypothesis testing and the estimation of spectral density. Section \ref{sec_5} presents empirical results based on selected models. To demonstrate the practical applicability of our method, we apply it to stock market data in Section \ref{application}, assessing the efficacy of our approach.

\section{Preliminaries}\label{pre}
\begin{definition}\label{Def4}
[Uniformly tight]\cite{billingsley1968convergence}
   Let $(X_n)$ be a sequence of $(S, \rho)$ valued random variables. $(X_n)$ is said to be \textbf{uniformly tight} if, for every $\epsilon > 0$, there exists a compact set $K_{\epsilon}\subseteq S$ such that $P(X_n \in K_{\epsilon}) > 1 - \epsilon$ for all $n$.
\end{definition}

\begin{lemma}\label{Lem3}
Suppose $X_n, X$ and $Y_n$'s are all $(S, \rho)$ valued random variables.
    If $(X_n)$ converges weakly (in distribution) to $X$ and $\rho(X_n,Y_n)\xrightarrow{P}0$, then $(Y_n)$ converges weakly (in distribution) to $Y$.
\end{lemma}

\begin{lemma}\label{Lem4}
[Prokhorov's Theorem]\cite{billingsley1968convergence}
Let $(X_n)$ be a sequence of $(S,\rho)$-valued random variable, where $S$ is polish, defined on a probability space $(\Omega, \mathcal{F}, P)$. Then $(X_n)$ contains a subsequence that converges weakly if and only if $(X_n)$ is uniformly tight.
\end{lemma}

\begin{proposition}\label{prop1} 
\cite{stout1974almost}
    Let $\{X_i, i \geq 1\}$ be the coordinate representation of 
a stationary sequence. There exists a measure-preserving transformation $\phi$
on $(\R_{\infty}, C_{\infty}, P)$ such that $X_1(\omega) = X_1(\omega)$, $X_2(\omega) = X(\phi(\omega))$, $\ldots$, $X_n(\omega) = X(\phi^n(\omega)), \ldots$ for all $\omega\in \R_{\infty}$. 
\end{proposition}

\begin{proposition}\label{prop3}
    [The pointwise ergodic theorem for stationary sequences]  \cite{stout1974almost}
    Let $\{X_i, i \geq 1\}$ be stationary with $\E| X_1| < \infty$. Then  $\sum_{i=1}^n X_i/n \rightarrow E[X_1| I]$ \ a.s. If in addition ${X_i , i \geq 1}$ is ergodic, $\sum_{i=1}^n X_i/n \rightarrow E[X_1]$ \ a.s.
\end{proposition}

\begin{proposition}\label{prop2} 
\cite{stout1974almost}
    Let $\{X_i, i \geq 1\}$ be stationary ergodic and $ \phi: R_{\infty} \rightarrow R$  be measurable. Let
$Y_i = \phi(X_i, X_{i+1}, \ldots), \forall i\geq 1$. Then ${Y_i, i \geq 1}$ is stationary ergodic.
\end{proposition}

\begin{lemma}\label{lemma_4}
    \cite{stout1974almost} 
    Let $\{X_i, i \geq 1\}$ be independent identically distributed. Then $\{X_i, i \geq 1\}$ is stationary ergodic.
\end{lemma}

\section{Theoretical results}\label{theory}
In this section, we establish the multivariate version of Theorem 21.1 of  \cite{billingsley1968convergence} for the case when $X_t = g(\xi(t), \xi(t-1), \ldots)$ is a d-dimensional sequence of random variables, with $g$ being measurable and $\epsilon$'s being m-dependent sequence of random variables, that satisfies the assumptions in Theorem 21.1 of \cite{billingsley1968convergence}. The following theorem helps us to perform hypothesis testing under the assumption that $H_0$ is true. Specifically, we show that the statistic $\sup_{t\in [0,1]}\sum_{i=1}^d (W_i^0(t))^2,$ where $ W_i^0(t)$'s are standard brownian-bridges in $[0,1]$, is the statistic which will be used for hypothesis testing and details of this statistics has been discussed in \cite{Kiefer}, where the determination of its quantile values have also been discussed. We have used those values for obtaining the empirical results in Section \ref{sec_5}.   

\begin{theorem}\label{thm1}
Consider a sequence $\{\xi(t)\}_{t \in \mathbb{Z}}$ of identically distributed, $m$-dependent, $d$- dimensional random variables with $\mathbb{E}[||\xi(t)||^2] < \infty$. Let $\{X_t\}_{t \in \mathbb{Z}}$ be a strictly stationary sequence of $d$-dimensional random variables defined as $X_t = g(\xi(t), \xi(t-1), \ldots)$, where $g$ is a measurable function. Denote $X_{tl} = g(\xi(t), \ldots, \xi(t-l))$. Assume $\mathbb{E}[\xi(t)] = 0$ and $\mathbb{E}[||X_t||^2] < \infty$ for all $t \in \mathbb{Z}$. Assume further that $$\sum_{l\geq 1}\left(\mathbb{E}\left[||X_0-X_{0l}||^2\right]\right)^{1/2}<\infty.$$ Then, 
$$\Sigma = \sum_{j\in \Z}Cov(X_0, X_j)$$ converges absolutely coordinatewise and if $\Sigma$ is positive-definite, then 
$$\frac{1}{\sqrt{N}}\sum_{i=1}^{[Nt]}X_i\xrightarrow{D[0,1]}W_{\Sigma}(t).$$

    \begin{proof}
        Let us denote
        \begin{align*}
            S_N(t) = \frac{1}{\sqrt{N}}\sum_{i=1}^{[Nt]}X_i.
        \end{align*}
 For a fixed $V\neq 0$, we denote $Y_i=V^TX_i=V^Tg(\xi(i), \xi({i-1}), \ldots)$ and $Y_{il}=V^TX_{il}=V^Tg(\xi(i),\ldots, \xi({i-l}))$, where $X_{il}=g(\xi(i),\ldots, \xi({i-l}))$.
Next, we have by the assumption above, that
\begin{align*}
    \sum_{l\geq 1}\left(\mathbb{E}\left[(Y_i-Y_{il})^2\right]\right)^{1/2} &= \sum_{l\geq 1}\left(\mathbb{E}\left[(V^T(X_i-X_{il}))^2\right]\right)^{1/2}\\
    &\leq ||V||_{\infty}\sum_{l\geq 1}\left(\mathbb{E}\left[||X_i-X_{il}||^2\right]\right)^{1/2} < \infty
\end{align*}
and $\mathbb{E}[V^TX_i]= 0.$
Hence $Y_i$ satisfies the assumptions in Theorem 21.1 in \cite{billingsley1968convergence}. Hence for $0\leq s \leq t\leq 1$, \begin{align*}
    S_N(s,t) = \frac{1}{\sqrt{N}}\sum_{i=[Ns]+1}^{[Nt]}Y_i\xrightarrow{D}N(0,\sigma_Y^2(t-s)),
\end{align*}
where
$\sigma_Y^2 = \sum_{j\in \Z} \mathbb{E}[Y_0Y_j] = \sum_{j\in \Z} V^TCov(X_0, X_j)V.$ (We have assumed that the matrix inside is positive definite, which implies that $\sigma_Y(t-s)>0$ for $V \neq 0$.) Since we can choose $V$ arbitrarily and by using similar arguments, we will have the same result. Hence by Cramer-Wald device, we have

$$\frac{1}{\sqrt{N}}\sum_{i=[Ns]+1}^{[Nt]}X_i\xrightarrow{D}W_{\Sigma} (t-s),$$ where $\Sigma = \sum_{j\in \Z}Cov(X_0, X_j)$ and $W_{\Sigma}$ is Wiener process with $\Sigma$ being the covariance matrix. For a fixed $0< s < t$, we consider the r.v s 
$$\Hat{S}_N(s,t) = \frac{1}{\sqrt{N}}\left(\sum_{i = [Ns] + 1}^{[Nt]}X_{i\left([Ns] + 1 -i\right)}\right)$$
and
$$S_{N,m}(s) = \frac{1}{\sqrt{N}}\left(\sum_{i = 1}^{[Ns]-m}X_{i}\right)$$ for some large $N$.
Therefore, $\hat{S}_N(s,t)$ and $\hat{S}_{N,m}(s)$ are independent random variables.
Now, upon using Markov inequality and Minkoswki inequality, we get
\begin{align*}
    \mathbb{P}\left(||\hat{S}_N(s,t) - S_{N}(s,t)|| > \delta\right) &\leq \frac{1}{\delta^2}\mathbb{E}\left[||\hat{S}_N(s,t) - S_{N}(s,t)||^2\right]\\
    & \leq \frac{1}{N\delta^2} \left( \sum_{i=[Ns]+1}^{[Nt]}\left(\mathbb{E}\left[||X_{i\left([Ns]+1-i\right)}-X_i||^2\right]\right)^{1/2} \right)^2\\
    & \leq \frac{1}{N\delta^2} \left( \sum_{i=[Ns]+1}^{[Nt]}\left(\mathbb{E}\left[||X_{1\left([Ns]+1-i\right)}-X_1||^2\right]\right)^{1/2} \right)^2\\
    & \leq \frac{1}{N\delta^2} \left( \sum_{l\geq0}\left(\mathbb{E}\left[||X_{1l}-X_1||^2\right]\right)^{1/2} \right)^2 \xrightarrow{N\xrightarrow{}\infty} 0,
\end{align*}
\begin{align*}
    \mathbb{P}\left(||S_N(t) - S_{N,m}(t)||>\delta\right) &\leq \frac{1}{N\delta^2}\mathbb{E}\left[||S_N(t) - S_{N,m}(t)||^2\right]\\
    &\leq \frac{1}{N\delta^2}\mathbb{E}\left[\left|\left|\sum_{i=[Nt]-m+1}^{[Nt]}X_i\right|\right|^2\right]\\
    & = \frac{m^2}{N\delta^2}\mathbb{E}\left(||X_i||^2\right)\xrightarrow{N\xrightarrow{}\infty}0.
\end{align*}
Hence, we have $||S_N(s) - S_{N.m}(t)||\xrightarrow{P}0$ and 
$||\hat{S}_N(s,t) - S_{N}(s,t)||\xrightarrow{P}0$. then, lemma \ref{Lem3} implies that
$$S_{N,m}(s) = \frac{1}{\sqrt{N}}\sum_{i=1}^{[Ns]}X_i\xrightarrow{D}N(0,\Sigma s)$$ and finally by Continuous Mapping Theorem, it follows that
\begin{align*}
    \left(\frac{1}{\sqrt{N}}\sum_{i=1}^{[Ns]}X_i, \frac{1}{\sqrt{N}}\sum_{i=1}^{[Nt]}X_i \right)\xrightarrow{D}\left(W_{\Sigma}(s), W_{\Sigma}(t)\right).
\end{align*}
Similarly, we can shows that the convergence in distribution takes place for any finite dimensional $(t_1,\ldots,t_k) \in [0,1]$, and tightness will follow when we choose $V = e_i$'s where $e_i$ is the vector with $1$ in the $i$-th position and zeros in all other places. This yields that the coordinates form a uniformly tight sequence of random variables and by Lemma \ref{Lem4}, it finally gives the uniform tightness of the $\{X_i\}_{i\geq 1}$. Hence, we have 
$$\frac{1}{\sqrt{N}}\sum_{i=1}^{[Nt]}X_i\xrightarrow{D[0,1]}W_{\Sigma}(t).$$
\end{proof}  
\end{theorem}
Then by Theorem \ref{thm1}, it follows that the CUSUM statistics converges weakly in Skorohod topology to $W^0_{\Sigma}(t)$, that is
\begin{align*}
  \frac{1}{\sqrt{N}}\left(\sum_{i=1}^{[Nt]}X_i-t\sum_{i=1}^{N}X_i\right)\xrightarrow{D[0,1]}W^0_{\Sigma}(t),  
\end{align*}
where $W^0_{\Sigma}(t)$ is Brownian bridge with $Cov(W^0_{\Sigma}(s), W^0_{\Sigma}(t)) = (\min{\{s,t\}} - st)\Sigma.$
Hence for a fixed $t$, we have $$(\Tilde{S}_N(t))^T\hat{\Sigma}^{-1}_N(\Tilde{S}_N(t))\xrightarrow{D[0,1]} \sum_{i=1}^{d}(W_i^0(t))^2,$$ where $$\Tilde{S}_N(t) = \frac{1}{\sqrt{N}}\left(\sum_{i=1}^{[Nt]}X_i-t\sum_{i=1}^{N}X_i\right)$$ and $W_i^0(t)$'s are standard one-dimensional Brownian-bridges in $[0,1]$ and $\hat{\Sigma}_N\xrightarrow{P}\Sigma$, and we will find such estimators in Theorem \ref{thm3}. 

Let $f\in D[0,1]$ and consider 
$N_{\epsilon}(f) = \{g\in D[0,1]: d_0(f,g)<\epsilon\}$ and for $\forall g\in N_{\epsilon}(f)$. By Definition \ref{Def4}, there exists a continuous strictly increasing function $\lambda:[0,1]\xrightarrow{}[0,1]$ such that $\lambda(0)=0$ and $\lambda(1)=1$ and $\sup_{t\in[0,1]} |f(t)-g(\lambda (t))|<3/2\epsilon$ and this implies that $|\sup_{t\in[0,1]} f(t) - \sup_{t\in[0,1]} g(t)|\leq 3/2\epsilon$. Hence, the function $\sup_{t\in[0,1]}(.)$ is a continuous function on $D[0,1]$. Thus, applying Continuous Mapping Theorem, we have 
$$\sup_{t\in[0,1]}(\Tilde{S}_N(t))^T\hat{\Sigma}^{-1}_N(\Tilde{S}_N(t))\xrightarrow{D} \sup_{t\in[0,1]}\sum_{i=1}^{d}(W_i^0(t))^2.$$

Here, we consider the linear process defined by $$X_t = \sum_{k \geq 0}C_k \xi(t-k),$$ where $\{\xi(t)\}_{t \in \mathbb{Z}}$ represents $m$-dependent identically distributed random variables. To satisfy the conditions stated in Theorem \ref{thm1}, we need

\begin{align*}\label{eqn_rev1}
     \sum_{l\geq 1}\mathbb{E} \left\{\left|\left|\sum_{m>l}C_m\xi({-m})\right|\right|^2\right\}^{1/2} < \infty.\\
\end{align*}

Again, using Minkowski Inequality, we have 
\begin{align*}
     \sum_{l\geq 1}\mathbb{E} \left\{\left|\left|\sum_{m>l}C_m\xi({-m})\right|\right|^2\right\}^{1/2} \leq \sum_{l\geq 1}\sum_{m>l}\left|\left|C_m\right|\right|\mathbb{E} \left\{\left|\left|\xi({-m})\right|\right|^2\right\}^{1/2}.\\
\end{align*}

Given that $\{\xi(t)\}_{-\infty < t < \infty}$ are identically distributed and $\mathbb{E} \left\{||\xi(0)||^2\right\}^{1/2} < \infty$, if
\begin{align*}
\sum_{l\geq 1}\sum_{m>l}||C_m|| < \infty,
\end{align*}
then the condition in Theorem \ref{thm1}, namely, $\sum_{l\geq 1}\mathbb{E}\left\{||X_0 - X_{0l}||^2\right\}^{1/2} < \infty$, is satisfied. This condition can be further extended to a condition based on the entries of the matrices $C_i$'s in the following manner:
\begin{align*}
     \sum_{l\geq 1}\sum_{m>l}\left|\left|C_m\right|\right| &\leq \sum_{l\geq 1}\sum_{m>l}\left(\sum_{i=1}^{d}\sum_{j=1}^{d}\left|(C_m)_{ij}\right|^2\right)^{1/2} \\
     & \leq \sum_{l\geq 1}\sum_{m>l}\left(\sum_{i=1}^{d}\sum_{j=1}^{d}\left|(C_m)_{ij}\right|\right) \\
     & \leq \sum_{i=1}^d\sum_{j=1}^{d}\left(\sum_{l\geq 1}\sum_{m>l}\left|(C_m)_{ij}\right|\right). \\
\end{align*}
Consequently, if we assume that the sum of the entries of the matrices $C_i$ converges absolutely, i.e., specifically $\sum_{l\geq 1}\sum_{m>l}\left|(C_m)_{ij}\right| < \infty$ for all $1 \leq i,j \leq d$, then the aforementioned condition gets satisfied. Thus, it is evident that under $H_0$, the $X_i$'s constitute a strictly stationary sequence of random variables, as required in Theorem \ref{thm1}. In the next theorem, we prove that under $H_a$, that is, when there exists one change-point and certain assumptions, the power of the test is asymptotically $1$ and we give a consistent estimator for the change-point.  
\begin{theorem}\label{thm2}
Suppose $T^* = [k^*N]$, for some $k^* \in (0,1)$. Let $\{X_t\}_{t \in \mathbb{Z}}$ be a strictly stationary and ergodic sequence with $\mathbb{E}[X_0]=0$ and $\mathbb{E}[||X_0||^2]<\infty$. Also, consider another strictly stationary and ergodic sequence $\{X^*_t\}_{t \in \mathbb{Z}}$ with $\mathbb{E}[X^*_0] \neq \mathbb{E}[X_0]$ and $\mathbb{E}[||X^*_0||^2]<\infty$. Define $Y_t = X_t$ for all $t \in [0, T^*]$ and $Y_t = X^*_t$, for all $t \in [T^*, N]$. In other words, $T^*$ represents the change-point in the mean of the data $\{Y_t\}_{t \in \mathbb{Z}}$. Then, both $$\argmax_{t\in[0,1]}\left|\left|\Tilde{S}_{N}(t)\right|\right|$$ and  $$\argmax_{t\in[0,1]}\frac{1}{n}\Tilde{S}_{N}(t)\Hat{\Sigma}_N^{-1}\Tilde{S}_{N}(t)$$ are consistent estimators of the point $k^*\in(0,1)$, where 
$\Tilde{S}_N(t) = \frac{1}{\sqrt{N}}\left(\sum_{i=1}^{[Nt]}Y_i-t\sum_{i=1}^{N}Y_i\right)$.
\begin{remark}
    The above theorem provides us a consistent estimator of $k^*$ and hence the estimate for the change-point will be $[k^*N]$, where $N$ is the length of the data.
\end{remark}
\begin{remark}\label{remark2}
    The assumption of ergodicity in the above theorem is quite reasonable from the fact that $X_t$'s and $X^*_t$'s are both of the form \ref{mod_1}, which is a measurable function of $\xi(t)$'s. And, if the observations come from the model \ref{mod_1} which we have used to simulate results, $\xi(t)$'s are again measurable functions of $Z(t)$'s which are independent and identically distributed and hence stationary ergodic by Lemma \ref{lemma_4} and hence it follows by Proposition \ref{prop2}. 
\end{remark}
    \begin{proof}
    By Proposition \ref{prop2}, we have  for a strictly stationary and ergodic sequence,
        \begin{align*}
            \sum_{i=1}^N \frac{X_i}{N}\xrightarrow{a.s} \mathbb{E}[X_1]. 
        \end{align*}
        This implies, by using the definiton  of almost sure convergence that for a fixed $T_0\in[0, k^*/2]$ and given $\delta>0$, $\hat{\epsilon} = k^*\epsilon>0$, there exists $N_0\in \mathbb{N}$ such that
        \begin{align*}
            & \mathbb{P}\left(\left|\left| \frac{1}{NT_0}\sum_{k=1}^{[NT_0]}X_i \right|\right|<k^*\epsilon, \forall N\geq [T_0N_0] \right) \geq 1-\delta\\
            \implies& \mathbb{P}\left(\left|\left| \frac{1}{NT}\sum_{k=1}^{[NT]}X_i \right|\right|<k^*\epsilon, \forall n\geq n_0, \forall T\in[T_0, k^*] \right) \geq 1-\delta\\
            \implies& \mathbb{P}\left(\left|\left| \frac{1}{N}\sum_{k=1}^{[NT]}X_i \right|\right|<\hat{\epsilon}, \forall N\geq N_0, \forall T\in[T_0, k^*] \right) \geq 1-\delta\\.
        \end{align*}
        In particular, we have 
        \begin{align*}       \mathbb{P}\left(\left|\left| \frac{1}{N}\sum_{k=1}^{[Nk^*]}X_i \right|\right|<\hat{\epsilon}, \forall N\geq N_0 \right) \geq 1-\delta.
        \end{align*}
        Hence using the previous observations, we get that 
        \begin{align*}
        &\mathbb{P}\left(\left|\left| \frac{1}{N}\sum_{k=[NT] + 1}^{[Nk^*]}X_i \right|\right|<2\hat{\epsilon}, \forall N\geq N_0, \forall T\in[T_0, k^*] \right) \geq 1-\delta  \\
          \implies &\mathbb{P}\left(\left|\left| \frac{1}{N}\sum_{k=[NT] + 1}^{[Nk^*]}X_i \right|\right|<2\hat{\epsilon}, \forall N\geq N_0, \forall T\in[T_0, k^*] \right) \geq 1-\delta  \\
          \implies &\mathbb{P}\left(\left|\left| \frac{1}{N}\sum_{k= 1}^{[Nk^*] - [NT]}X_i \right|\right|<2\hat{\epsilon}, \forall n\geq N_0, \forall T\in[T_0, k^*] \right) \geq 1-\delta  \\
          \implies &\mathbb{P}\left(\left|\left| \frac{1}{N}\sum_{k= 1}^{[NT^{'}]} 
          X_i \right|\right|<2\hat{\epsilon}, \forall N\geq N_0, \forall T^{'}\in[0, k^*/2] \right) \geq 1-\delta \\
          \implies &\mathbb{P}\left(\left|\left| \frac{1}{N}\sum_{k= 1}^{[NT]}X_i \right|\right|<2\hat{\epsilon}, \forall N\geq N_0, \forall T\in[0, k^*] \right) \geq 1-\delta\\
          \implies &\lim_{{N \to \infty}} \mathbb{P} \left( \sup_{{T \in [0,k^*]}}  \left| \left| \frac{1}{N} \sum_{{k=1}}^{{\lfloor NT \rfloor}} X_i \right|\right| < 2\hat{\epsilon}   \right) = 1
        \end{align*}
and similarly by the assumptions that there is another strictly stationary and ergodic sequence of random variables 
$\{X_i^*\}_{i\geq 0}$ with $\mathbb{E}[X_i^*] = \mathbb{E}[X_0^*], \forall i \in \mathbb{Z}$ and $X_i^* = Y_i, \forall i\geq [Nk^*]$.
\begin{align*}
\lim_{N\xrightarrow{}\infty}\mathbb{P}\left(\sup_{T\in[0,k^*]}\left|\left| \frac{1}{N}\sum_{k= [Nk^*]+ 1}^{[Nt]}X^*_i - (t-k^*)\mathbb{E}[X_0^*]\right|\right|>2\hat{\epsilon} \right) = 0.
\end{align*}
We thus obtain $$\mathbb{P}\left(\sup_{t\in[0,1]}\left|\left|\frac{1}{\sqrt{N}}\Tilde{S}_N(t)-g(t)\right|\right|> \epsilon\right)\xrightarrow{N\xrightarrow{}\infty}0,$$
where 
$$\Tilde{S}_N(t) = \frac{1}{\sqrt{N}}\left(\sum_{i=1}^{[Nt]}Y_i-t\sum_{i=1}^{N}Y_i\right)$$ and
$$g(t) = \begin{cases} 
      t(1-k^*)(\mathbb{E}[X_0] - \mathbb{E}[X^*_0]), & \text{if } t \in [0, k^*] \\
      k^*(1-t)(\mathbb{E}[X_0] - \mathbb{E}[X^*_0]), & \text{if } t \in [k^*, 1]. 
   \end{cases}$$
Therefore, $\sup_{t\in[0,1]}\left|\frac{1}{\sqrt{N}}\Tilde{S}_N(t)-g(t)\right|\xrightarrow[P]{N\xrightarrow{}\infty}0$.
Hence, for a given $\epsilon>0$, $\exists \text{ } \delta >0, N_0\in\mathbb{N}$ such that 
\begin{equation*}
\mathbb{P}\left(\sup_{t\in[0,1]}\left|\left|\frac{1}{\sqrt{N}}\Tilde{S}_N(t)-g(t)\right|\right|< \delta\right)>1-\epsilon,\forall N\geq N_0.    
\end{equation*}
Hence, $\forall N\geq N_0$ with probability at least $1-\epsilon$, 
\begin{equation*}
   \sqrt{(k^*)^2-2\delta} < \argmax_{t\in[0,k^*]}\left|\left|\frac{1}{\sqrt{N}}\Tilde{S}_n(t)\right|\right|^2\leq k^* 
\end{equation*}
and similarly,
\begin{equation*}
    k^*\leq \argmax_{t\in[0,k^*]}\left|\left|\frac{1}{\sqrt{N}}\Tilde{S}_N(t)\right|\right|^2<\sqrt{(k^*)^2+2\delta}.
\end{equation*}
Thus, for $N\geq N_0$, we have an injective function $\Hat{\Delta}(\delta)$ dependent only on $k^*$ such that 
\begin{align*}
    \mathbb{P}\left(\left| \argmax_{t\in[0,k^*]}\left|\left|\frac{1}{\sqrt{N}}\Tilde{S}_N(t)\right|\right|^2-k^*\right|<\Hat{\Delta}(\delta)\right) > 1-\epsilon,
\end{align*}
and hence $$\argmax_{t\in[0,1]}\left|\left|\Tilde{S}_{N}(t)\right|\right|$$ and similarly $$\argmax_{t\in[0,1]}\frac{1}{N}\Tilde{S}_{N}(t)\Hat{\Sigma}_N^{-1}\Tilde{S}_{N}(t)$$ are consistent estimators of $k^*\in(0,1)$.
    \end{proof}
\end{theorem}
\begin{remark}
Theorem \ref{thm2} can be extended to the case of multiple change-points. Say, there are $k$ change-points in the data, rather than just a single change-point. Using the same proof under the assumption that there are $k$ such stationary-ergodic sequences, we can show that
$$\frac{1}{N}\Tilde{S}_{N}(t)\Hat{\Sigma}_n^{-1}\Tilde{S}_{N}(t)$$ has local etremas at the change-points.    
\end{remark}
Now, the final task is to find a consistent estimator of $\Sigma$ in Theorem \ref{thm1} which will be used for hypothesis testing.
Suppose the periodogram of $\{X_1,\ldots, X_N\}$, where ${X_t}$ is a stationary time series with $\mathbb{E}[X_t] = \mu$ and covariance matrices $\mathbb{E}[X_{t+h}X^T_t]-\mu\mu^T = \Gamma(h)$ has absolutely summable components. Under this condition, ${X_t}$ has a continuous spectral density given by
$$f(\lambda) = \frac{1}{2\pi}\sum_{h=-\infty}^{\infty} \Gamma(h)e^{-ih\lambda}, \quad -\pi \leq \lambda \leq \pi.$$

For $\omega_j=\frac{2\pi j}{N}$, where $-[(N-1)/2]\leq j \leq [N/2]$ the discrete Fourier transform of  $\{X_1, \ldots, X_N\}$ is defined by 
\begin{align*}
    W(\omega_j) = \frac{1}{\sqrt{N}}\sum_{n=1}^{N}X_n e^{in\omega_j}
\end{align*}
and let $I_N(\omega_j)$ denote the periodogram of $\{X_1, \ldots, X_N\}$, defined by $I_N(\omega_j) = W(\omega_j)W(\omega_j)^*$.
Thus, $$I_{N}(\omega_j) = \frac{1}{N}\left(\sum_{n=1}^{N}X_n e^{in\omega_j}\right)\left(\sum_{n=1}^{N}X_n e^{in\omega_j}\right)^*.$$
The definition of $I_N(\omega)$ is then extended to $\omega\in [-\pi, \pi]$ by the following relation
\begin{align*}
    I_N(\omega) &= 
    \begin{cases}
        I_N(g(N, \omega)), & \text{if } \omega \geq 0, \\
        I_N^*(g(N, -\omega)), & \text{if } \omega < 0,
    \end{cases}
\end{align*}
where $g(N, \omega)$, where $\omega \in [0,\pi]$, is the nearest multiple of $2\pi/N$ to $\omega$.

In the following theorem, we estimate the spectral density $f(\omega)$,
where $\omega \in [0,\pi]$. We consider an estimator for $f(\omega)$, as $\Hat{f}(\omega) = (2\pi)^{-1}\sum_{|k|\leq h_N} K_N(k)I_N(g(N, \omega)+\omega_{k}),$ where $\{h_N\}$ is a sequence of positive integers and $\{K_N(.)\}$ is a sequence of weight functions. We also impose the following conditions on $\{h_N\}$ and $\{K_N(.)\}$ :
\begin{align}\label{weight}
    &h_N\xrightarrow{}\infty \text{  and  } h_N/\sqrt{N}\xrightarrow{} 0 \text{ as } N\xrightarrow{} \infty,\nonumber\\
    &K_N(k) = K_N(-k), \hspace{0.5in} K_N(k)\geq 0, \text{ for all } k \nonumber\\
    & \sum_{|k|\leq h_N} K_N(k) = 1\nonumber\\
     \text{ and }& \sum_{|k|\leq h_N} K_N^2(k) \xrightarrow{N\xrightarrow{}\infty} 0 .
\end{align}
\begin{remark}
 The above assumptions are satisfied by the Simple moving average kernel
    \begin{align}
        K_N(k) = \frac{1}{(2h_N + 1)}.\mathbb{I}_{\{|k|\leq h_N\}},
    \end{align}
    where $h_N$ is the chosen sequence of bandwidth. 
\end{remark}

\begin{theorem}\label{thm3}
Under the assumptions on $X_t$ based on Theorem \ref{thm1} and the weight functions $K_N(k)$ or under $H_0$, $\hat{f}(\omega) = (2\pi)^{-1}\sum_{|k|\leq h_N} K_N(k)I_N(g(N, \omega)+\omega_{k})$ is a consistent estimator of $f(\omega)$, for $-\pi \leq \omega \leq \pi$. Consequently, there exists a consistent estimator $\hat{\Sigma}_N$ for $\Sigma$. 
\begin{proof}
To prove this, we follow similar steps as in Section 10.3 in \cite{brockwell_davis} where they prove the univariate version for linear processes with independent innovations.
We similarly define the discrete Fourier transform of  $\{\xi(1), \ldots, \xi(N)\}$ as 
\begin{align*}
    W_{N,\xi}(\lambda) = \frac{1}{\sqrt{N}}\sum_{n=1}^{N}\xi(n) e^{in\lambda}.
\end{align*}
For $\omega_j = \frac{2\pi j}{N}$, where $-[(N-1)/2]\leq j \leq [N/2]$, 
$\left(I_{N,\xi}(\omega_j)\right)$ denotes the periodogram of $\{\xi(1), \ldots, \xi(N)\}$, defined by $I_{N,\xi}(\omega_j) = W_{N,\xi}(\omega_j)W_{N,\xi}(\omega_j)^*$.
Thus, 
$$\left(I_{N,\xi}(\omega_j)\right)_{pq} = \frac{1}{N}\left(\sum_{n=1}^{N}\xi_{p}(n) e^{in\omega_j}\right)^*\left(\sum_{n=1}^{N}\xi_{q}(n) e^{in\omega_j}\right).$$
\begin{align*}
\mathbb{E}\left[ \left(I_{N,\xi}(\omega_j)\right)_{pq}  \left(I_{N,\xi}(\omega_k)\right)_{rs}\right] &= \frac{1}{N^2}\sum_{s,t,u,v=1}^{N}\mathbb{E}\left[ \xi_p(s)\xi_q(t)\xi_r(u)\xi_s(v)e^{i(t-s)\omega_j}e^{i(v-u)\omega_k}\right] 
\end{align*}
The last sum can be decomposed into four types of sums (since we have assumed that $\mathbb{E}\left[\xi(t)\right] = 0$) as follows:
\begin{align*}
\sum\sum\sum\sum \mathbb{E}\left[ \xi_p(s)\xi_q(t)\xi_r(u)\xi_s(v)\right]e^{i(t-s)\omega_j}e^{i(v-u)\omega_k} \hspace{0.9in}\text{(i)} \\
+ \sum\sum\sum\sum \mathbb{E}\left[ \xi_p(s)\xi_q(t)]\E[\xi_r(u)\xi_s(v)\right]e^{i(t-s)\omega_j}e^{i(v-u)\omega_k} \hspace{0.9in}\text{(ii)}\\
+ \sum\sum\sum\sum \mathbb{E}\left[ \xi_p(s)\xi_q(u)]\E[\xi_r(t)\xi_s(v)\right]e^{i(t-s)\omega_j}e^{i(v-u)\omega_k} \hspace{0.9in}\text{(iii)}\\
+\sum\sum\sum\sum \mathbb{E}\left[ \xi_p(s)\xi_q(v)]\E[\xi_r(t)\xi_s(v)\right]e^{i(t-s)\omega_j}e^{i(v-u)\omega_k} \hspace{0.9in}\text{(iv)}.
\end{align*}
For the type-i sum, since the $\xi(i)$'s are m-dependent random variables, it should be of $O(N)$. For the type-iii sum, we consider the following expressions,
\begin{align*}
    &\sum_{t = m+1}^{N-m}\sum_{\alpha=-m}^{m} \mathbb{E}\left[\xi_p(t)\xi_q(t+\alpha)\right]e^{i\omega_jt}e^{i\omega_k(t+\alpha)} \\
    &=\sum_{t = m+1}^{N-m}\sum_{\alpha=-m}^{m} \mathbb{E}\left[\xi_p(t)\xi_q(t+\alpha)\right]e^{i(\omega_j+\omega_k)t}e^{i\omega_k\alpha}\\
    &=\sum_{t = 1}^{N}\sum_{\alpha=-m}^{m} \mathbb{E}\left[\xi_p(t)\xi_q(t+\alpha)\right]e^{i(\omega_j+\omega_k)t}e^{i\omega_k\alpha} \\
    &\hspace{2.5in} - \sum_{t = 1}^{m}\sum_{\alpha=-m}^{m} \mathbb{E}\left[\xi_p(t)\xi_q(t+\alpha)\right]e^{i(\omega_j+\omega_k)t}e^{i\omega_k\alpha} \\
    &\hspace{2.5in} - \sum_{t = N-m+1}^{N}\sum_{\alpha=-m}^{m} \mathbb{E}\left[\xi_p(t)\xi_q(t+\alpha)\right]e^{i(\omega_j+\omega_k)t}e^{i\omega_k\alpha}\\
    &=\left(\sum_{t = 1}^{N} e^{i(\omega_j+\omega_k)t}\right)\left(\sum_{\alpha=-m}^{m}\mathbb{E}\left[\xi_p(t)\xi_q(t+\alpha)\right]e^{i\omega_k\alpha}\right)\\
    &\hspace{2.5in} - \sum_{t = 1}^{m}\sum_{\alpha=-m}^{m} \mathbb{E}\left[\xi_p(t)\xi_q(t+\alpha)\right]e^{i(\omega_j+\omega_k)t}e^{i\omega_k\alpha} \\
    &\hspace{2.5in} - \sum_{t = N-m+1}^{N}\sum_{\alpha=-m}^{m} \mathbb{E}\left[\xi_p(t)\xi_q(t+\alpha)\right]e^{i(\omega_j+\omega_k)t}e^{i\omega_k\alpha}\\
    &=\left(\sum_{t = 1}^{N} e^{i(\omega_j+\omega_k)t}\right)\left(\sum_{\alpha=-m}^{m}\mathbb{E}\left[\xi_p(t)\xi_q(t+\alpha)\right]e^{i\omega_k\alpha}\right) + O(1),
\end{align*}
In the third equality of the above expression, we have used the fact that $\xi(n)'s$ is a stationary sequence of random variables and so we can decompose it into products which doesn't depend on $t$. And the remaining part of the expression in the third equality if of $O(1)$ since there are $m$ terms, where $m$ is merely a constant and hence again by stationary property of $\xi(n)'s$, we have the rest to be of constant order.
This means the following:
\begin{align*}
\text{If } \omega_j = \omega_k = 0 \text{ or } \omega_j = \omega_k = \pi, & \text{ then }  = O(N) ;\\
\text{Otherwise,} & \text{ then }  = O(1).
\end{align*}
Therefore, the following is true for type-iii sum: 
\begin{align*}
\text{If } \omega_j = \omega_k = 0 \text{ or } \omega_j = \omega_k = \pi, & \text{ then }  = O(N^2) ;\\
\text{Otherwise,} & \text{ then }  = O(1).
\end{align*}
Similar arguments will work for the type-iv sum as well and will give the same result. We thus have 
\begin{align*}
    &\text{Cov}\left( (I_{N,\xi}(\omega_j))_{pq}, (I_{N,\xi}(\omega_k))_{rs}  \right) \\&=\mathbb{E}\left( \left(I_{N,\xi}(\omega_j)\right)_{pq}\left(I_{N,\xi}(\omega_k)\right)_{rs}\right) - \mathbb{E} \left(I_{N,\xi}(\omega_j)\right)_{pq}\mathbb{E} \left(I_{N,\xi}(\omega_j)\right)_{pq} \\
    &= \mathbb{E}\left( \left(I_{N,\xi}(\omega_j)\right)_{pq}\left(I_{N,\xi}(\omega_k)\right)_{rs}\right)\\& \hspace{1.5in}- \frac{1}{N^2}\sum_{s,t,u,v=1}^{N} \mathbb{E}\left[ \xi_p(s)\xi_q(t)]\E[\xi_r(u)\xi_s(v)\right]e^{i(t-s)\omega_j}e^{i(v-u)\omega_k}.
\end{align*}

Hence combining the above orders, we get
\begin{align}\label{eqn_th3_2}
\text{Cov}\left( (I_{N,\xi}(\omega_j))_{pq}, (I_{N,\xi}(\omega_k))_{rs}  \right) &= 
\begin{cases}
O\left(\frac{1}{N}\right), & \text{if } 0<\omega_j \neq \omega_k<\pi, \\
O(1), & \text{if } 0 \leq \omega_j = \omega_k \leq  \pi.
\end{cases}
\end{align}

Likewise, the discrete Fourier transform of  $\{X_1, \ldots, X_N\}$ is defined as 
\begin{align*}
    &W_{N,X}(\lambda) = \frac{1}{\sqrt{N}}\sum_{n=1}^{N}X_n e^{in\lambda}\\
    &= \frac{1}{\sqrt{N}}\sum_{n=1}^{N}e^{in\lambda}\left\{\sum_{k\geq0}C_k \xi(t-k)\right\}\\
    &= \sum_{k\geq0}C_k e^{ik\lambda}\frac{1}{\sqrt{N}}\sum_{n=1}^{N}\xi(n)e^{in\lambda} + \frac{1}{\sqrt{N}}\sum_{k\geq0}C_k e^{ik\lambda}R_{k,N}(\lambda),
\end{align*}
where 
\begin{align*}
    \left(R_{k,N}(\lambda)\right)_p = \mathbb{E}\left[-\sum_{n=1}^{N}\xi_p(n)e^{in\lambda} - \sum_{n=-k+1}^{N-k}\xi_p(n)e^{in\lambda} \right].
\end{align*}
Let $\left(I_{N,X}(\lambda)\right)$ denote the periodogram of $\{X_1, \ldots, X_N\}$, defined by \\ $I_{N,X}(\lambda) = W_{N,X}(\lambda)W_{N,X}(\lambda)^*$. 

Thus 
\begin{align}\label{eq_th3_1}
I_{N,X}(\lambda) &= \frac{1}{N}\left(\sum_{n=1}^{N}X_{n} e^{in\lambda}\right)\left(\sum_{n=1}^{N}X_{n} e^{in\lambda}\right)^*\nonumber\\
& = \left(\sum_{k\geq0}C_k e^{ik\lambda}\frac{1}{\sqrt{N}}\sum_{n=1}^{N}\xi(n)e^{in\lambda} + \frac{1}{\sqrt{N}}\sum_{k\geq0}C_k e^{ik\lambda}R_{k,N}(\lambda)\right)\nonumber\\
&\hspace{1in}\times \left(\sum_{k\geq0}C_k e^{ik\lambda}\frac{1}{\sqrt{N}}\sum_{n=1}^{N}\xi(n)e^{in\lambda} + \frac{1}{\sqrt{N}}\sum_{k\geq0}C_k e^{ik\lambda}R_{k,N}(\lambda)\right)^*\nonumber\\
&=\left(\sum_{k\geq0}C_k e^{ik\lambda}\right)I_{N,\xi}(\lambda)\left(\sum_{k\geq0}C_k e^{ik\lambda}\right)^*\nonumber\\
&\hspace{1in} + \left(\frac{1}{\sqrt{N}}\sum_{k\geq 0} C_ke^{ik\lambda}R_{k,N}(\lambda)\right)\left(\frac{1}{\sqrt{N}}\sum_{k\geq 0} C_ke^{ik\lambda}R_{k,N}(\lambda)\right)^* \nonumber\\
&\hspace{1in}+ \left(\frac{1}{\sqrt{N}}\sum_{k\geq 0} C_ke^{ik\lambda}R_{k,N}(\lambda)\right)\left(\sum_{k\geq0}C_k e^{ik\lambda}\frac{1}{\sqrt{N}}\sum_{n=1}^{N}\xi(n)e^{in\lambda}\right)^*\nonumber\\
&\hspace{1in} + \left(\sum_{k\geq0}C_k e^{ik\lambda}\frac{1}{\sqrt{N}}\sum_{n=1}^{N}\xi(n)e^{in\lambda}\right)\left(\frac{1}{\sqrt{N}}\sum_{k\geq 0} C_ke^{ik\lambda}R_{k,N}(\lambda)\right)^* \nonumber\\
&= \left(\sum_{k\geq0}C_k e^{ik\lambda}\right)I_{N,\xi}(\lambda)\left(\sum_{k\geq0}C_k e^{ik\lambda}\right)^* + R_N(\lambda), 
\end{align}
where 
\begin{align*}
   R_N(\lambda) &= \left(\frac{1}{\sqrt{N}}\sum_{k\geq 0} C_ke^{ik\lambda}R_{k,N}(\lambda)\right)\left(\frac{1}{\sqrt{N}}\sum_{k\geq 0} C_ke^{ik\lambda}R_{k,N}(\lambda)\right)^* \\
&\hspace{1in}+ \left(\frac{1}{\sqrt{N}}\sum_{k\geq 0} C_ke^{ik\lambda}R_{k,N}(\lambda)\right)\left(\sum_{k\geq0}C_k e^{ik\lambda}\frac{1}{\sqrt{N}}\sum_{n=1}^{N}\xi(n)e^{in\lambda}\right)^*\\
&\hspace{1in} + \left(\sum_{k\geq0}C_k e^{ik\lambda}\frac{1}{\sqrt{N}}\sum_{n=1}^{N}\xi(n)e^{in\lambda}\right)\left(\frac{1}{\sqrt{N}}\sum_{k\geq 0} C_ke^{ik\lambda}R_{k,N}(\lambda)\right)^*.
\end{align*}
By using Minkowski’s inequality,

\begin{equation}
\begin{split}
\max_{\lambda}\frac{1}{\sqrt{2\pi N}} \left(\mathbb{E}\left[\left(\left(\sum_{j=-\infty}^{\infty} C_j e^{ij\lambda} R_{j,N}(\lambda)\right)^* \left(\sum_{j=-\infty}^{\infty} C_j e^{ij\lambda} R_{j,N}(\lambda)\right)\right)^2\right]\right)^{1/4} \\
\leq \frac{1}{\sqrt{2\pi N}} \max_{\lambda} \sum_{j=-\infty}^{\infty} \left(\mathbb{E} \left[ \left|\left| C_j e^{ij\lambda} R_{j,N}(\lambda) \right|\right|^4 \right]\right)^{1/4} \\
\leq \frac{1}{\sqrt{2\pi N}} \max_{\lambda} \sum_{j=-\infty}^{\infty} \left|\left|C_j\right|\right| \left(\mathbb{E} \left[ \left|\left| R_{j,N}(\lambda) \right|\right|^4 \right]\right)^{1/4}
\end{split}
\end{equation}
and
\begin{equation}
\mathbb{E}\left[ \left|\left| R_{k,N}(\lambda) \right|\right|^4 \right] = \mathbb{E}\left[ \left( \sum_{m,n\in I} \xi_p\left(m\right)\xi_p\left(n\right)e^{i(m-n)\lambda}\right)^2 \right] = O\left(\min(|k|^2,|N|^2)\right)
\end{equation}

\text{where} $I = \{1,\ldots,N\} \Delta \{-k+1,\ldots,N-k\}$. Therefore by the assumption that
\begin{align*}
\sum_{l\geq 1}\sum_{m>l}||C_m|| < \infty,
\end{align*}
we have
\begin{align}\label{rev_eqn1}
    &\max_{\lambda}{\mathbb{E}\left[\left(\left(\frac{1}{\sqrt{2\pi N}}\sum_{j=-\infty}^{\infty} C_je^{ij\lambda}R_{j,N}(\lambda)\right)^*\left(\frac{1}{\sqrt{2\pi N}}\sum_{j=-\infty}^{\infty} C_je^{ij\lambda}R_{j,N}(\lambda)\right)\right)^{2}\right]}\nonumber\\&
    \leq O\left(\frac{1}{N^2}\right).
\end{align}
and
\begin{align*}
 \mathbb{E}\left[ \left| C(e^{i\omega_k})W_{N,\xi}(\omega_k)r_n^*(\omega_k) \right|^2_{ij} \right] &= \mathbb{E}\left[ \left| C(e^{i\omega_k})W_{N,\xi}(\omega_k)\right|^2_{i} \left|r_N^*(\omega_k) \right|^2_{j} \right] \\ 
 &\leq \mathbb{E}\left[ \left| C(e^{i\omega_k})W_{N,\xi}(\omega_k)\right|^4_{i} \right]^{1/2} \mathbb{E}\left[\left|r_N(\omega_k) \right|^4_{j} \right]^{1/2}\\
  &\leq \left|\left|C(e^{i\omega_k})\right|\right|^2 \mathbb{E}\left[ \left|W_{N,\xi}(\omega_k)\right|^4_{i} \right]^{1/2} \mathbb{E}\left[\left|r_N(\omega_k) \right|^4_{j} \right]^{1/2},
\end{align*}
where $r_N(\omega_k) = \frac{1}{\sqrt{N}}\sum_{j=-\infty}^{\infty} C_je^{ij\omega_k}R_{j,N}(\omega_k)$ for $\omega_k = \frac{2\pi k}{N}\in [0, \pi]$. As we already have seen that
\begin{align*}
&\text{Cov}\left( (I_N(\omega_k))_{pq}, (I_N(\omega_k))_{rs}  \right)\\ 
&=\mathbb{E}\left[ \left( (W_{N,\xi}(\omega_k))_p\overline{(W_{N,\xi}(\omega_k))_{q}} \right)\left( (W_{N,\xi}(\omega_k))_r\overline{(W_{N,\xi}(\omega_k))_s} \right) \right]\\
& \hspace{1in}- \mathbb{E}\left[ \left( (W_{N,\xi}(\omega_k))_p\overline{(W_{N,\xi}(\omega_k))_q} \right)\right]\mathbb{E}\left[\left(( W_{N,\xi}(\omega_k))_r\overline{(W_{N,\xi}(\omega_k))_s} \right) \right]\\
&=O\left(\frac{1}{N}\right) + O(1)
\end{align*}
and
\begin{align*}
 \mathbb{E}\left[ \left( (W_{N,\xi}(\omega_k))_p\overline{(W_{N,\xi}(\omega_k))_q} \right)\right]\mathbb{E}\left[\left(( W_{N,\xi}(\omega_k))_r\overline{(W_{N,\xi}(\omega_k))_s} \right) \right]= O(1)  
\end{align*}
since the sequence 
$\{\xi\left(n\right)\}_{-\infty < n < \infty}$ are $m$-dependent random variables, we have
\begin{align*}
    \mathbb{E}\left[ \left( (W_{N,\xi}(\omega_k))_p\overline{(W_{N,\xi}(\omega_k))_q} \right)\right] = \frac{1}{2\pi N}\sum_{j, k =1}^{N} \mathbb{E}\left[\xi_p\left(j\right)\xi_q\left(k\right)\right]e^{i(j-k)\omega_k} = O\left(1 \right).
\end{align*}
Hence,
\begin{align*}
\mathbb{E}\left[ \left( W_p(\omega_k)\overline{W_q(\omega_k)} \right) \left( W_r(\omega_k)\overline{W_s(\omega_k)} \right) \right] &= O\left(\frac{1}{N}\right) + O(1), 
\end{align*}
which implies that
\begin{align*}
\max_{\omega_k\in[0,\pi]}\mathbb{E}\left[ \left|W(\omega_k)\right|^4_{i} \right] = O\left(\frac{1}{N}\right) + O(1),
\end{align*}
\begin{align}\label{eqn_thm3_3}
\max_{\omega_k\in[0,\pi]} \mathbb{E}\left[ \left| C(e^{i\omega_k})W(\omega_k)r_N^*(\omega_k) \right|^2_{ij} \right] &\leq \left|\left|C(e^{i\omega_k})\right|\right|^2 \left( O\left(\frac{1}{N}\right) + O\left(1 \right) \right)^{1/2} \left(O\left(\frac{1}{N^2}\right) \right)^{1/2}\nonumber
\\ & = O\left(\frac{1}{N}\right)
\end{align}
and
\begin{align}\label{eqn_thm3_5}
&\max_{\lambda}\mathbb{E}\left|\left(\frac{1}{\sqrt{N}}\sum_{k\geq 0} C_ke^{ik\lambda}R_{k,N}(\lambda)\right)\left(\frac{1}{\sqrt{N}}\sum_{k\geq 0} C_ke^{ik\lambda}R_{k,N}(\lambda)\right)^*\right|^2_{ij} \nonumber\\
&=\max_{\lambda} \mathbb{E}\left|\left(\frac{1}{\sqrt{N}}\sum_{k\geq 0} C_ke^{ik\lambda}R_{k,N}(\lambda)\right)_i\overline{\left(\frac{1}{\sqrt{N}}\sum_{k\geq 0} C_ke^{ik\lambda}R_{k,N}(\lambda)\right)_j}\right|^2 \nonumber\\
&=\max_{\lambda}\mathbb{E}\left|\left(\frac{1}{\sqrt{N}}\sum_{k\geq 0} C_ke^{ik\lambda}R_{k,N}(\lambda)\right)_i\right|^2\left|\left(\frac{1}{\sqrt{N}}\sum_{k\geq 0} C_ke^{ik\lambda}R_{k,N}(\lambda)\right)_j\right|^2 \nonumber\\
&\leq \max_{\lambda}\mathbb{E}\left|\left(\frac{1}{\sqrt{N}}\sum_{k\geq 0} C_ke^{ik\lambda}R_{k,N}(\lambda)\right)\right|^2\left|\left(\frac{1}{\sqrt{N}}\sum_{k\geq 0} C_ke^{ik\lambda}R_{k,N}(\lambda)\right)\right|^2\nonumber\\
&\leq\max_{\lambda} \left(\mathbb{E}\left|\left(\frac{1}{\sqrt{N}}\sum_{k\geq 0} C_ke^{ik\lambda}R_{k,N}(\lambda)\right)\right|^4\mathbb{E}\left|\left(\frac{1}{\sqrt{N}}\sum_{k\geq 0} C_ke^{ik\lambda}R_{k,N}(\lambda)\right)\right|^4\right)^{1/2} \nonumber\\
&= O\left(\frac{1}{N^2}\right).
\end{align} 
This is true from the observation obtained in (\ref{rev_eqn1}).
Finally by Minkowski's inequality, we have 
\begin{align*}
&\max_{\omega_k\in[0,\pi]}
\left(\mathbb{E}\left[ \left|R_N(\omega_k) \right|_{ij}^2 \right]\right)^{1/2} \\&
\leq \max_{\omega_k\in[0,\pi]}\mathbb{E}\left(\left[\left(\frac{1}{\sqrt{N}}\sum_{k\geq 0} C_ke^{ij\omega_k}R_{j,N}(\omega_k)\right)\left(\frac{1}{\sqrt{N}}\sum_{k\geq 0} C_ke^{ij\omega_k}R_{j,N}(\omega_k)\right)^*\right]_{ij}^2 \right)^{1/2}\\
&+ \max_{\omega_k\in[0,\pi]}\left(\mathbb{E}\left[\left(\sum_{k\geq 0} C_ke^{ij\omega_k}R_{j,N}(\omega_k)\right)\left(\sum_{k\geq0}C_k e^{ik\omega_k}\frac{1}{\sqrt{N}}\sum_{n=1}^{N}\xi(n)e^{in\omega_k}\right)^*\right]_{ij}^2\right)^{1/2}\\
&+ \max_{\omega_k\in[0,\pi]}\left(\mathbb{E}\left[\left(\sum_{k\geq0}C_k e^{ik\omega_k}\frac{1}{\sqrt{N}}\sum_{n=1}^{N}\xi(n)e^{in\omega_k}\right)\left(\sum_{k\geq 0} C_ke^{ij\omega_k}R_{j,N}(\omega_k)\right)^*\right]_{ij}^2\right)^{1/2}\\&\leq O\left(\frac{1}{\sqrt{N}}\right).
\end{align*}
Hence,
\begin{align}\label{eqn_thm_3}
  \max_{\omega_k\in[0,\pi]}\mathbb{E}\left[ \left|R_N(\omega_k) \right|_{ij}^2 \right]\leq O\left(\frac{1}{N}\right).  
\end{align}
Therefore from (\ref{eq_th3_1}), we have 
\begin{align*}
    & (I_{N,X}(\omega_k))_{ij} = \left(C(e^{i\omega_k})I_{N,Z}(\omega_k)C(e^{i\omega_k})^*\right)_{ij} + \left(R_N(\omega_k)\right)_{ij}\\
    & \hspace{0.85in} = \sum_{1\leq m,n\leq d} \left[\left(C(e^{i\omega_k})\right)_{im}\left(I_{N,Z}(\omega_k)\right)_{mn}\left(C(e^{i\omega_k})\right)_{nj}\right] + (R_N(\omega_k))_{ij}.
\end{align*}
which implies that
\begin{align*}
   &\text{Cov}\left((I_{N, X}(\omega_j))_{pq},
    (I_{N, X}(\omega_k))_{rs}\right) 
   \\&= \text{Cov}\bigg(\sum_{1\leq m,n\leq d} \left(C(e^{i\omega_j})\right)_{pm}\left(I_{N,\xi}(\omega_j)\right)_{mn}\left(C(e^{i\omega_j})\right)_{nq},\\& \hspace{1in} \sum_{1\leq m,n\leq d} \left(C(e^{i\omega_k})\right)_{rm}\left(I_{N,\xi}(\omega_k)\right)_{mn}\left(C(e^{i\omega_k})\right)_{ns}\bigg)\\& \hspace{1in} + \text{Cov}\left(\sum_{1\leq m,n\leq d} \left(C(e^{i\omega_j})\right)_{pm}\left(I_{N,\xi}(\omega_j)\right)_{mn}\left(C(e^{i\omega_j})\right)_{nq}, (R_N(\omega_k))_{rs}\right)\\& \hspace{1in}+ \text{Cov}\left((R_N(\omega_j))_{pq}, \sum_{1\leq m,n\leq d} \left(C(e^{i\omega_k})\right)_{rm}\left(I_{N,\xi}(\omega_k)\right)_{mn}\left(C(e^{i\omega_k})\right)_{ns}\right) \\& \hspace{1in}+ \text{Cov}\left((R_N(\omega_j))_{pq}, (R_N(\omega_k))_{rs}\right).
\end{align*}
Hence from the relations in  (\ref{eqn_th3_2}), (\ref{rev_eqn1}), (\ref{eqn_thm3_3}), (\ref{eqn_thm3_5}), (\ref{eqn_thm_3}) and using $\text{Cov}(\text{ . , . })$ function and Cauchy-Schawtz inequality, we obtain

\begin{align*}
    \text{Cov}\left((I_{N, X}(\omega_j))_{pq}, (I_{N, X}(\omega_k))_{rs}\right) &= 
    \begin{cases}
        O\left(\frac{1}{\sqrt{N}}\right) + O(1), & \text{if } 0 \leq \omega_j = \omega_k \leq \pi\\
        O\left(\frac{1}{\sqrt{N}}\right), & \text{if } \text{otherwise}.
    \end{cases}
\end{align*}
By the restriction assumed on $h_N$ based on (\ref{weight}) that $$\frac{h_N}{\sqrt{N}}\xrightarrow{N\xrightarrow{}\infty}0,$$ we have
\begin{align*}
    \max_{|k|\leq h_N}\left| g(N,\omega) + \omega_k - \omega \right| \xrightarrow {N\hookrightarrow\infty} 0.
\end{align*}
This implies that, by continuity of $f_{ij}(.)$ on the compact set $[0,\pi]$, and hence by uniform continuity, we have
\begin{align*}
   &\max_{|k|\leq h_N}\left| f_{ij}\left(g(N,\omega) + \omega_k \right) - f_{ij}\left(\omega\right) \right| \xrightarrow {N\hookrightarrow\infty} 0. 
\end{align*}
As
\begin{align*}
\bigl|\mathbb{E}\hat{f}_{ij}(\omega) - f_{ij}(\omega)\bigr| = \biggl|&\sum_{|k|\leq h_N} K_N(k)\bigl[(2\pi)^{-1} \mathbb{E}\left(I_{N,X}\left(g(N,\omega) + \omega_k\right)_{ij}\right)\\
& - f_{ij}\left(g(N,\omega) + \omega_k\right) + f_{ij}\left(g(N,\omega) + \omega_k \right)-f_{ij}\left(\omega \right)\bigr]\biggr|
\end{align*}
 by Proposition 10.3.1 in \cite{brockwell_davis}, for stationary sequence $X_t$'s, we have $$\max_{|k|\leq h_N}\bigl|(2\pi)^{-1} \mathbb{E}\left(I_{N,X}\left(g(N,\omega) + \omega_k\right)_{ij}\right)
 - f_{ij}\left(g(N,\omega) + \omega_k\right) \bigr|\xrightarrow{N\xrightarrow{}\infty}0.$$
Therefore, we have $\mathbb{E}\hat{f}(\omega)\xrightarrow{N\xrightarrow{}\infty}f(\omega)$. Hence, for $\omega \in (0, \pi)$, 
\begin{align*}
    &\text{Var}\left(\hat{f}_{pq}(\omega)\right) = (2\pi)^{-2}\sum_{|j|\leq h_N} K_N^2(j)\left((2\pi)^2\text{Var}\left(\left(I_{N,X}(g(N,\omega)+\omega_j)\right)_{pq}\right)\right)\\ &+ (2\pi)^{-2}\sum_{|j|\leq h_N}\sum_{|k|\leq h_N, k\neq j} K_N(j)K_N(k)\\&\hspace{0.5in}\times\bigg((2\pi)^2\text{Cov}\left(\left(I_{N,X}(g(N,\omega)+\omega_j)\right)_{pq},
    0\left( I_{N,X}(g(N,\omega)+\omega_k)\right)_{pq}\right)\bigg)\\
    &=\left(\sum_{|j|\leq h_N}K_N^2(j)\right)O\left(1\right) + (2h_N+1)\left(\sum_{|j|\leq h_N}K_N^2(j)\right)O\left(\frac{1}{\sqrt{N}}\right)\\
    &= \left(\sum_{|j|\leq h_N}K_N^2(j)\right)O\left(1\right) + o\left(\sum_{|j|\leq h_N}K_N^2(j)\right)
\end{align*}
This implies that
$$\text{Var}\left(\hat{f}_{pq}(\omega)\right)\xrightarrow{N\xrightarrow{}\infty}0 , \text{ when } 0<\omega<\pi.$$
For $\omega = 0$, we have
\begin{align*}
    \Hat{f}(0) &= (2\pi)^{-1}\sum_{|k|\leq h_N} K_N(k)I_N(\omega_{k}) \\&= (2\pi)^{-1}\sum_{0\leq k\leq h_N} K_N(k)I_N(\omega_{k}) + (2\pi)^{-1}\sum_{-h_N \leq k< 0} K_N(k)I_N(\omega_{k})
\end{align*} and lets denote the sums as $\hat{f_1}(0)$ and $\hat{f_2}(0)$, respectively. Hence,

\begin{align*}
 &\text{Var}\left(\hat{f_1}_{pq}(0)\right) = (2\pi)^{-2}\sum_{0\leq j\leq h_N} K_N^2(j)\left((2\pi)^2\text{Var}\left(\left(I_{N,X}(\omega_j)\right)_{pq}\right)\right)\\ &+ (2\pi)^{-2}\sum_{0\leq j\leq h_N}\sum_{0\leq k\leq h_N, k\neq j} K_N(j)K_N(k)\left((2\pi)^2\text{Cov}\left(\left(I_{N,X}(\omega_j)\right)_{pq},\left( I_{N,X}(\omega_k)\right)_{pq}\right)\right)\\
    &\leq \left(\sum_{|j|\leq h_N}K_N^2(j)\right)O\left(1\right) + (2h_N+1)\left(\sum_{|j|\leq h_N}K_N^2(j)\right)O\left(\frac{1}{\sqrt{N}}\right)\\
    &= \left(\sum_{|j|\leq h_N}K_N^2(j)\right)O\left(1\right) + o\left(\sum_{|j|\leq h_N}K_N^2(j)\right) \xrightarrow{N\xrightarrow{} \infty} 0   
\end{align*}
and similarly, $\text{Var}\left(\hat{f_2}_{pq}(0)\right)\xrightarrow{N\xrightarrow{} \infty}0$. Since $\hat{f}_{pq}(0)= \hat{f_1}_{pq}(0) + \hat{f_2}_{pq}(0)$, 
$$\text{Var}(\hat{f}_{pq}(0)) \leq \text{Var}(\hat{f_1}_{pq}(0)) + \text{Var}(\hat{f_2}_{pq}(0)) + 2 \sqrt{\text{Var}(\hat{f_1}_{pq}(0))\text{Var}(\hat{f_2}_{pq}(0))}$$ by Cauchy-Schwarz inequality,
Therefore,
$$\text{Var}\hat{f}_{pq}(0)\xrightarrow{N\xrightarrow{}\infty}0,$$
and hence, $$\mathbb{E}\left|\Hat{f}_{pq}(0) - f_{pq}(0)\right|^2 = \text{Var}\Hat{f}_{pq}(0) + \left|\mathbb{E}\Hat{f}_{pq}(0) - f_{pq}(0)\right|^2\xrightarrow{N\xrightarrow{}\infty}0,$$ as required.
This implies that $\Hat{f(0)}$ is a consistent estimator of $f(0).$ Hence, proved.
\end{proof}
\end{theorem}
\section{Algorithm and Empirical results}\label{sec_5}
In this section, $T$ represents the length of a given data. The empirical results were obtained by utilizing the critical values as detailed in the paper by \cite{Kiefer} and using Simple moving average kernel with bandwidth $h_T = T^{1/4}$, which satisfies the assumptions \ref{weight}. 
For the task of single change-point detection, the proposed method was applied to datasets spanning different time lengths, namely, $8000$, and $16000$, across a spectrum of $m$-dependence levels including $10$ and $20$. As expected, given the consistent nature of our estimator and its convergence to the requisite statistic discussed in \cite{Kiefer} for large $T$, we observe an improvement in the performance of the method with increasing $T$. Moreover, since the assumed constant-order dependence of $\xi(n)'s$ suggests a relationship, where $O(m)<< O(T^{1/4})$ as obtained in the proof of Theorem \ref{thm3}, the selection of $m$ becomes crucial in optimizing the method's efficacy relative to the dataset size $T$. Additionally, it is worth noting that the performance of the method may be influenced by the location of the actual change-point, given the utilization of the convergence of the average of $X_t's$ as $T$ tends to infinity in the proof of Theorem \ref{thm2}. Another potential factor impacting method efficacy is the magnitude of the change itself, which warrants consideration in our evaluation.

For generating the data which follows the assumptions in Theorem \ref{thm1}, we generated a sequence of independent and identically distributed random variables $\{Z(t)\}_{t\in \Z}$, and then generated $\xi(t) = \Tilde{f}(Z(t-m), \ldots, Z(t+m)), \forall$ $t \in \Z$, with some suitable function $\Tilde{f}(^.)$. Hence, the data we used for analysis are as follows:
\begin{align}\label{mod_1}
X_t = \sum_{k \geq 0}C_k \Tilde{f}(Z(t-k-m), \ldots,Z(t), \ldots, Z(t-k+m)),
\end{align} where $m$ is the dependence parameter of the $\xi(t)$'s. This is a specific model which generates such data and we have used this to obtain our simulation results.
\par
\begin{figure}[h]
  \centering
  \includegraphics[width=1\textwidth]{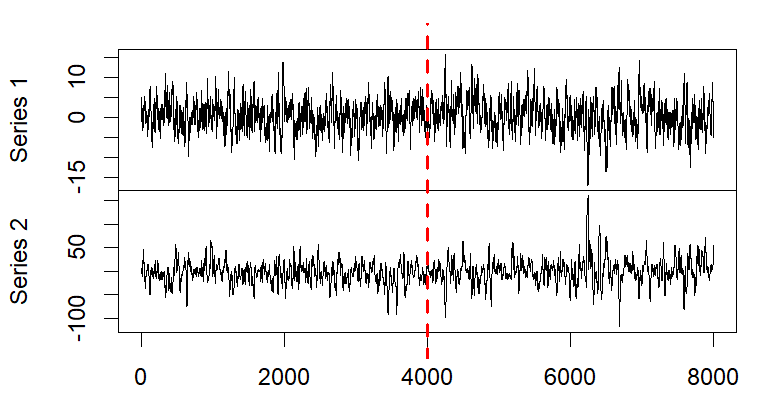}
  \caption{Data for the bi-variate with T = 8000 and change of mean being (0.2,0.5) at time T/2.}
  \label{plot_11}
\end{figure}

\begin{figure}[h]
  \centering
  \includegraphics[width=1\textwidth]{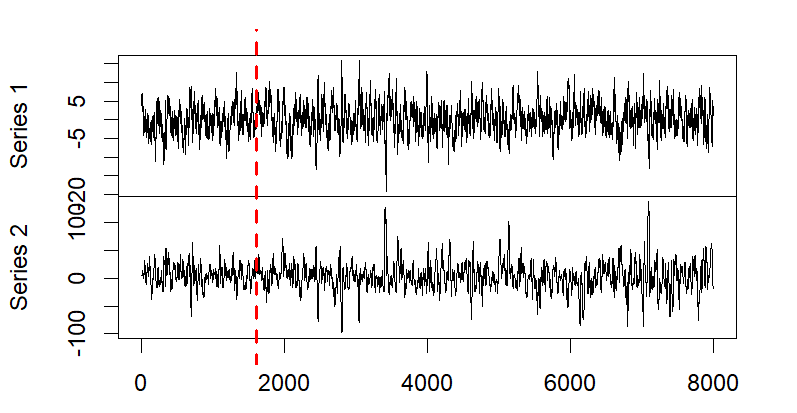}
  \caption{Data for the bi-variate with T = 8000 and change of mean being (0.2,0.5) at time T/5.}
  \label{plot_12}
\end{figure}

 For each combination in the tables, we have used 30 runs and took the average of the deviations of the estimates of the change-points from the actual change-point, absolute value of deviations and square of the deviations to estimate $\E\left[T^* - \Hat{T}\right]$ (deviation), $\E\left|T^*-\Hat{T}\right|$ (abs deviation) and $\E\left|T^*-\Hat{T}\right|^2$ (sqd deviation), where $\Hat{T}$ and $T^*$ are the estimated and actual change-points, respectively. While using the model in \ref{mod_1}, we used normal distribution with mean $\bold{0}$ and covariance matrix 
$${\begin{Bmatrix} 1&0.5 \\ 0.5&1 \end{Bmatrix}}$$
and
$${\begin{Bmatrix} 1&0.5&0.5&0.5&0.5 \\ 0.5&1&0.5&0.5&0.5 \\ 0.5&0.5&1&0.5&0.5 \\ 0.5&0.5&0.5&1&0.5 \\ 0.5&0.5&0.5&0.5&1 \end{Bmatrix}},$$ for generating $\{Z(t)\}_{t\in \Z}$ for the bi-variate and five-variate data cases, respectively, and observed the performance of the method based on two magnitude of change in the mean $(0.5, 0.2)$ and $(0.5, 1.2)$ for the bi-variate case and $(0.5, 0.2, 0.2, 0.5, 0.2)$ and $(0.5, 1.2, 0.5, 0.5, 0.5)$ for the five-variate case. The Figures \ref{plot_11} and \ref{plot_12} show the simulated bi-variate data with the actual change-point at $T/2$ and $T/5$, respectively. It is clear from the Figures that it is difficult to identify the change-points, mainly because of the spread of the data. 

For each simulation, we perform the hypothesis testing, based on the results obtained in Theorem \ref{thm1} at a significance level of $0.05$ and on an average, we couldn't reject the null hypothesis $1$ or $2$ out of $30$ simulations for each case. And once the null hypothesis gets rejected, we estimate the change-point based on the data simulated. We discuss some of the important observations from the Tables \ref{tab:data}-\ref{tab:data2_similar}. For change of (0.5, 0.2) with $m$-dependence of 20, we see that performance of the estimator, $\argmax_{t\in[0,1]}\frac{1}{T}\Tilde{S}_{T}(t)\Hat{\Sigma}_T^{-1}\Tilde{S}_{T}(t)$ in Theorem \ref{thm2} decreases when the actual change-point is at $T/5$; for example, the mean absolute deviation for 8000 length data increase to 1143 from 752.7, which is for the case when the location is $T/2$. But, when the change is increased to $(0.5, 1.2)$, there is not much of a change in the performance in terms of the mean absolute deviation and in both the cases, they perform well. This is mainly because of the increase in the magnitude of change in mean; the estimation is not affected to a large extent due to the change in location in the actual change-point. For the $m$-dependence 10, the performance of the estimator is good for all the cases, except for the case when the location of the change-point is $T/5$ and change in mean of $(0.5, 0.2)$, but the decline of the performance is not as bad as that of the case of $m$-dependence $20$. This is mainly due to the choice of $m$. Similarly, the estimation is good enough for the change in mean being $(0.5, 1.2)$ and the variation in performance is also low. We see similar observations for the datasets of length 16000, but the change in the performance is not that drastic for the $m$-dependence of 20 when the location of the change-point is changed from $T/2$ to $T/5$. We also present some more simulation results based on $m$-dependence of 30. For this case, we see that the performance of the estimator for the case when the change in mean is $(0.5,1.2)$ is still very good, especially for dataset of length 16000. So, we can expect to obtain an increase in the accuracy of the method by allowing larger dependence based on the chosen value of $T$ and the amount of change. Similar observations are seen in case of five-variate case when we alter the $m$-dependence, change in location of change-point and magnitude of change in the mean. 

\begin{figure}[h]
  \centering
  \includegraphics[width=1\textwidth]{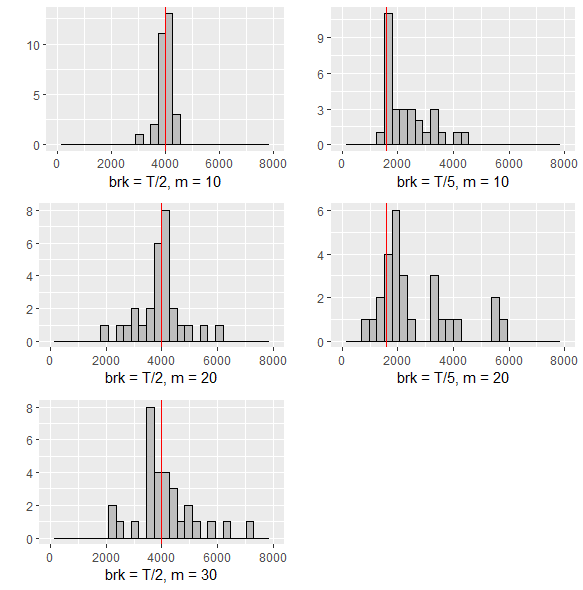}
  \caption{Conditional distribution of the change-point: T = 8000, Mean change = (0.2,0.5), x-axis represents the time index and y-axis represents the frequency of the estimates of the change-points.}
  \label{plot_1}
\end{figure}

\begin{figure}[h]
  \centering
  \includegraphics[width=1\textwidth]{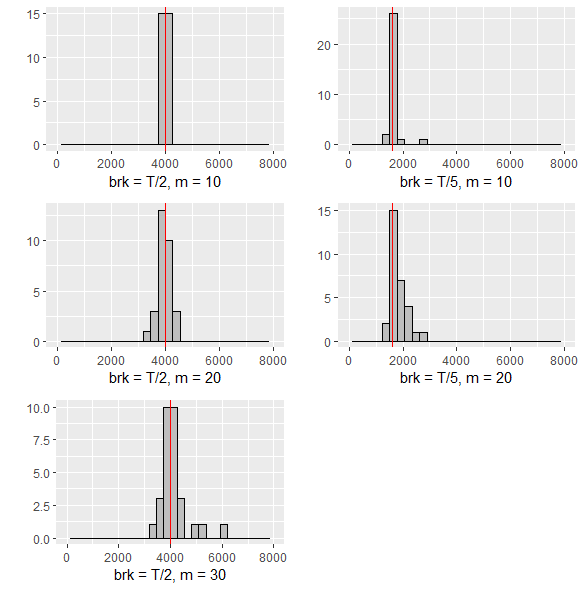}
  \caption{Conditional distribution of the change-point estimates: T = 8000, Mean change = (0.5,1.2), x-axis represents the time index and y-axis represents the frequency of the estimates of the change-points.}
  \label{plot_2}
\end{figure}

\begin{figure}[h]
  \centering
  \includegraphics[width=1\textwidth]{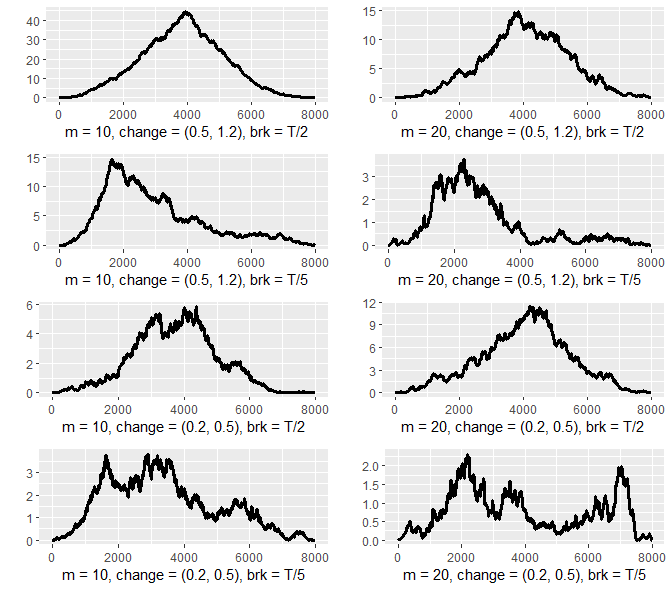}
  \caption{Performance of the estimator $\argmax_{t\in[0,1]}\frac{1}{T}\Tilde{S}_{T}(t)\Hat{\Sigma}_T^{-1}\Tilde{S}_{T}(t)$ based on one simulation for each of the cases with T = 8000.}
  \label{estimator}
\end{figure}

To show a graphical view of the results in Tables \ref{tab:data}-\ref{tab:data2_similar}, we present Figures \ref{plot_1} and \ref{plot_2} to show the distribution of the estimated change-points once the null hypothesis is rejected. We also show the Figure \ref{estimator} for the combinations in Table \ref{tab:data}. The plots clearly demonstrate the improvement in performance when the magnitude of change is increased and $m$-dependence is decreased. Similar plots can also be obtained for other tables as well. 

\begin{table}[htbp]\label{Table2_8000}
  \centering
  \caption{Performance evaluation of the estimator on the bi-variate data, T = 8000}
  \resizebox{1.0\textwidth}{!}{%
    \begin{tabular}{cccccc}
    \toprule
    change & m-dependence & location & deviation & abs deviation & sqd-deviation \\
    \midrule
    0.5, 0.2 & 10 & T/2 & 8.3 & 154.1 & 238.5 \\
            & 10 & T/5 & -714.7 & 728.9 & 813.02 \\
            & 20 & T/2 & 87.5 & 564.5 & 857.1 \\
            & 20 & T/5 & -1021.1 & 1143 & 1356.9 \\
            & 30 & T/2 & -85.33 & 752.7 & 1097 \\
    0.5, 1.2  & 10 & T/2 & -3.6 & 25.4 & 38.6 \\
            & 10 & T/5 & -74.5 & 99.1 & 216.3 \\
            & 20 & T/2 & 45.8 & 144.2 & 214.5 \\
            & 20 & T/5 & -226.3 & 280.8 & 358.6 \\
            & 30 & T/2 & -95.5 & 352.6 & 551.2 \\
    \bottomrule
    \end{tabular}%
  }
  \label{tab:data}%
\end{table}%

\begin{table}[htbp]\label{Table2_16000}
  \centering
  \caption{Performance evaluation of the estimator on the bi-variate data, T = 16000.}
  \resizebox{1.0\textwidth}{!}{%
    \begin{tabular}{cccccc}
    \toprule
    change & m-dependence & location & deviation & abs deviation & sqd-deviation \\
    \midrule
    0.5, 0.2 & 10 & T/2 & 6.1 & 99.6 & 178.2 \\
            & 10 & T/5 & -415 & 522.6 & 687 \\
            & 20 & T/2 & 44.3 & 315 & 459.3 \\
            & 20 & T/5 & -452 & 745 & 854 \\
            & 30 & T/2 & 96 & 311.2 & 425 \\
    0.5, 1.2  & 10 & T/2 & 7 & 12 & 29 \\
            & 10 & T/5 & -52 & 75 & 124 \\
            & 20 & T/2 & 18.3 & 67.8 & 110.2 \\
            & 20 & T/5 & -204 & 213 & 336 \\
            & 30 & T/2 & 55 & 157 & 248 \\
    \bottomrule
    \end{tabular}%
  }
  \label{tab:data2}%
\end{table}%

\begin{table}[htbp]\label{Table5_8000}
  \centering
  \caption{Performance evaluation of the estimator on the five-variate data, T = 8000.}
  \resizebox{1.05\textwidth}{!}{%
    \begin{tabular}{cccccc}
    \toprule
    change & m-dependence & location & deviation & abs deviation & sqd-deviation \\
    \midrule
    0.5, 0.2, 0.2, 0.5, 0.2 & 10 & T/2 & 11.2 & 179.1 & 275.6 \\
            & 10 & T/5 & -686.6 & 755.0 & 863.7 \\
            & 20 & T/2 & 93.6 & 573.4 & 848.8 \\
            & 20 & T/5 & -913.4 & 1054.5 & 1594.8 \\
            & 30 & T/2 & -63.4 & 662.7 & 902.0 \\
    0.5, 1.2, 0.5, 0.5, 0.5  & 10 & T/2 & -6.8 & 20.4 & 32.3 \\
            & 10 & T/5 & -69.8 & 107.1 & 196.3 \\
            & 20 & T/2 & 53.1 & 137.3 & 250.5 \\
            & 20 & T/5 & -241.8 & 312.8 & 375.4 \\
            & 30 & T/2 & -88.3 & 322.9 & 560.6 \\
    \bottomrule
    \end{tabular}%
  }
  \label{tab:data_similar_8000}%
\end{table}%

\begin{table}[htbp]\label{Table5_16000}
  \centering
  \caption{Performance evaluation of the estimator on the five-variate dataset, T = 16000.}
  \resizebox{1.05   \textwidth}{!}{%
    \begin{tabular}{cccccc}
    \toprule
    change & m-dependence & location & deviation & abs deviation & sqd-deviation \\
    \midrule
    0.5, 0.2, 0.2, 0.5, 0.2 & 10 & T/2 & 9.7 & 105.8 & 208.1 \\
            & 10 & T/5 & -408.2 & 568.7 & 679.2 \\
            & 20 & T/2 & 52.9 & 332.8 & 478.1 \\
            & 20 & T/5 & -405.4 & 713.2 & 874.9 \\
            & 30 & T/2 & 121.7 & 351.7 & 471.6 \\
    0.5, 1.2, 0.5, 0.5, 0.5  & 10 & T/2 & 6.7 & 16.2 & 27.3 \\
            & 10 & T/5 & -46.8 & 87.3 & 138.7 \\
            & 20 & T/2 & 26.9 & 77.9 & 115.3 \\
            & 20 & T/5 & -287.1 & 296.7 & 425.9 \\
            & 30 & T/2 & 61.2 & 146.9 & 255.5 \\
    \bottomrule
    \end{tabular}%
  }
  \label{tab:data2_similar}%
\end{table}%

\section{Application}\label{application}
\begin{figure}
    \centering
    \includegraphics[width=0.9\textwidth]{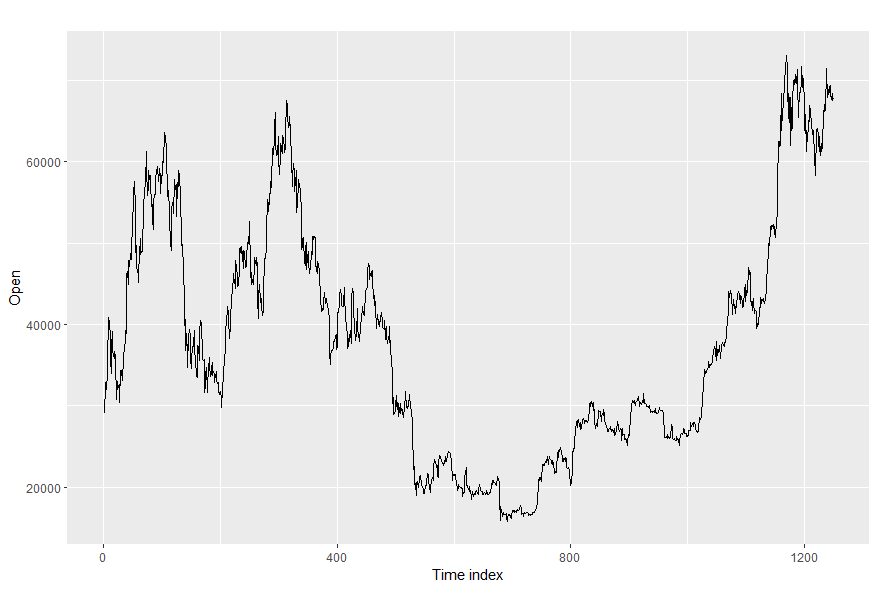}
    \caption{Daily opening price of bitcoin data from 1st January, 2021 to 31st May, 2024.}
    \label{bit-open}
\end{figure}

\begin{figure}
    \centering
    \includegraphics[width=0.9\textwidth]{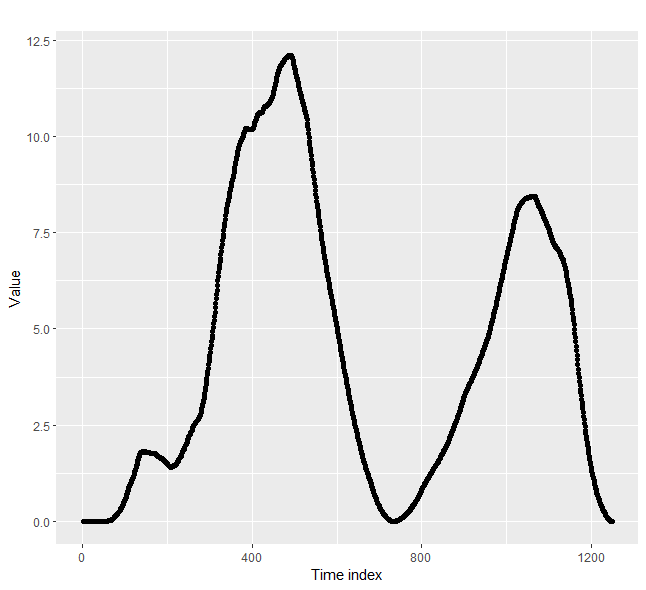}
    \caption{Graph of $\frac{1}{T}\Tilde{S}_{T}(t)\Hat{\Sigma}_T^{-1}\Tilde{S}_{T}(t)$ for the bit-coin data showing the local maximas and minimas, corresponding to the estimated multiple change-points.}
    \label{bit-est}
\end{figure}

We implemented our method on a comprehensive dataset of Bitcoin prices from \href{https://archive.ics.uci.edu/}{UCI Repository}, encompassing daily records from January 1, 2021, to May 31, 2024. This dataset included crucial financial metrics: Opening Price, Closing Price, High Price, Low Price, and Adjusted Price. Our analysis revealed a striking pattern, the local maxima and minima of $\frac{1}{T}\Tilde{S}_{T}(t)\Hat{\Sigma}_T^{-1}\Tilde{S}_{T}(t)$ in Figure \ref{bit-est} identified by our method coincided precisely with the most significant spikes and dips in the Bitcoin stock price in Figure \ref{bit-open}.
The graphical representation of our results  was nothing short of astonishing. We try to estimate the change-points of the bitcoin data by looking into its local maximas and minimas. Figure \ref{bit-est} vividly illustrated how our method pinpointed the moments of dramatic market shifts. Each peak and trough in our data was accurately matched with the local maxima and minima, showcasing the method's precision. These moments of sharp market movements, characterized by sudden surges and plummets, were captured with remarkable clarity, affirming the robustness and reliability of our approach.

This visual evidence underscores the efficacy of our method in tracking and predicting critical points in Bitcoin's volatile market, providing invaluable insights for traders and analysts seeking to navigate the ever-fluctuating landscape of cryptocurrency.

\section{Conclusion.}
This work has explored the theoretical results, where we first present Theorem \ref{thm1}. The assumptions made in this theorem are mainly the assumptions that are made to perform hypothesis testing on whether there is actually a shift in the mean or not. For this, we use the statistic in Theorem \ref{thm1}, to perform the hypothesis testing. It is also worth mentioning that to be able to calculate the statistic, we found a consistent estimator $\hat{\Sigma}_N$ for $\Sigma$ in Theorem \ref{thm3}. Once the testing is done and the null hypothesis gets rejected, we use Theorem \ref{thm2} to get the consistent estimator of the change-point under certain assumptions like ergodicity, mentioned in the theorem. These assumptions are justified by Remark \ref{remark2}, where we discuss why this assumption fits the with the model \ref{mod_1} with which we generate the empirical results. The results points out some observations such as the effect of increase in length of the data, value of $m$ in the $m$-dependence, location of the actual change-point and magnitude of the shift on the performance of the estimator. Although the method doesn't include any assumption on the dependence structure within the components, a potential limitation of the proposed method is that the \(\Sigma\) matrix mentioned in Theorem \ref{thm1} may not always be positive-definite. This might pose some issue while performing the hypothesis testing as well as estimation of the change-point.

Despite the progress made in this field, there are still problems that are of interest for further exploration and development. Future research could focus on refining existing methods to handle more complex data structures, such as high-dimensional datasets or those with non-linear dependencies. An immediate question which may arise is finding the change-points based on change in distribution of innovations which in turn causes change in distribution of the data. Mattenson \cite{matteson2014nonparametric} has explored a method to detect change-points based on change in distribution of data assuming that the observations are independent, but using the linear process model with change in distribution of the innovations might be an alternative way of looking into the problem involving dependence structure of data. 

Ultimately, the continued advancement of multivariate change-point detection methods holds great promise for improving our understanding of dynamic systems and in facilitating informed decision-making in a wide range of applications. Thus, this area of research provides great opportunities for researchers to push the boundaries of knowledge in this direction.

\section{Acknowledgements}
    I would like to express my sincere gratitude to Prof. N. Balakrishnan for his revision of this manuscript. I also extend my heartfelt thanks to my friend, 
 \href{https://ankankar-zargon.github.io/}{\textcolor{black}{Ankan Kar}} [ORCID: \href{https://orcid.org/0009-0003-7084-0961}{\includegraphics[scale=0.1]{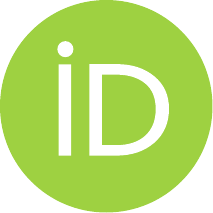}}], for his assistance with the manuscript writing and Latex Code. 
 
\bibliography{references}

\begin{thebibliography}{11}

\bibitem[\protect\citeauthoryear{Billingsley}{1968}]{billingsley1968convergence}
\begin{bbook}[author]
\bauthor{\bsnm{Billingsley},~\bfnm{Patrick}\binits{P.}}
(\byear{1968}).
\btitle{Convergence of Probability Measures}.
\bpublisher{John Wiley \& Sons, New York}.
\end{bbook}
\endbibitem

\bibitem[\protect\citeauthoryear{Brockwell and Davis}{1991}]{brockwell_davis}
\begin{bbook}[author]
\bauthor{\bsnm{Brockwell},~\bfnm{Peter~J.}\binits{P.~J.}} \AND \bauthor{\bsnm{Davis},~\bfnm{Richard~A.}\binits{R.~A.}}
(\byear{1991}).
\btitle{Time Series: Theory and Methods},
\bedition{2nd} ed.
\bpublisher{Springer, New York}.
\end{bbook}
\endbibitem

\bibitem[\protect\citeauthoryear{Gao, Yang and Yang}{2020}]{gao2020change}
\begin{barticle}[author]
\bauthor{\bsnm{Gao},~\bfnm{Wei}\binits{W.}}, \bauthor{\bsnm{Yang},~\bfnm{Haizhong}\binits{H.}} \AND \bauthor{\bsnm{Yang},~\bfnm{Lu}\binits{L.}}
(\byear{2020}).
\btitle{Change points detection and parameter estimation for multivariate time series}.
\bjournal{Soft Computing}
\bvolume{24}
\bpages{6395--6407}.
\end{barticle}
\endbibitem

\bibitem[\protect\citeauthoryear{Horv{\'a}th and Hu{\v{s}}kov{\'a}}{2012}]{horvath2012change}
\begin{barticle}[author]
\bauthor{\bsnm{Horv{\'a}th},~\bfnm{Lajos}\binits{L.}} \AND \bauthor{\bsnm{Hu{\v{s}}kov{\'a}},~\bfnm{Marie}\binits{M.}}
(\byear{2012}).
\btitle{Change-point detection in panel data}.
\bjournal{Journal of Time Series Analysis}
\bvolume{33}
\bpages{631--648}.
\end{barticle}
\endbibitem

\bibitem[\protect\citeauthoryear{Kiefer}{1959}]{Kiefer}
\begin{barticle}[author]
\bauthor{\bsnm{Kiefer},~\bfnm{J.}\binits{J.}}
(\byear{1959}).
\btitle{{K-Sample Analogues of the Kolmogorov-Smirnov and Cramer-V. Mises Tests}}.
\bjournal{The Annals of Mathematical Statistics}
\bvolume{30}
\bpages{420 -- 447}.
\end{barticle}
\endbibitem

\bibitem[\protect\citeauthoryear{Kuncheva}{2011}]{kuncheva2011change}
\begin{barticle}[author]
\bauthor{\bsnm{Kuncheva},~\bfnm{Ludmila~I}\binits{L.~I.}}
(\byear{2011}).
\btitle{Change detection in streaming multivariate data using likelihood detectors}.
\bjournal{IEEE Transactions on Knowledge and Data Engineering}
\bvolume{25}
\bpages{1175--1180}.
\end{barticle}
\endbibitem

\bibitem[\protect\citeauthoryear{Matteson and James}{2014}]{matteson2014nonparametric}
\begin{barticle}[author]
\bauthor{\bsnm{Matteson},~\bfnm{David~S}\binits{D.~S.}} \AND \bauthor{\bsnm{James},~\bfnm{Nicholas~A}\binits{N.~A.}}
(\byear{2014}).
\btitle{A nonparametric approach for multiple change point analysis of multivariate data}.
\bjournal{Journal of the American Statistical Association}
\bvolume{109}
\bpages{334--345}.
\end{barticle}
\endbibitem

\bibitem[\protect\citeauthoryear{Page}{1955}]{pagetest}
\begin{barticle}[author]
\bauthor{\bsnm{Page},~\bfnm{ES}\binits{E.}}
(\byear{1955}).
\btitle{A test for a change in a parameter occurring at an unknown point}.
\bjournal{Biometrika}
\bvolume{42}
\bpages{523--527}.
\end{barticle}
\endbibitem

\bibitem[\protect\citeauthoryear{S.~Kotz and Vidakovic}{1998}]{johnson1998encyclopedia}
\begin{bbook}[author]
\bauthor{\bsnm{S.~Kotz},~\bfnm{C.~Read}\binits{C.~R.} \bsuffix{N.~Balakrishnan}} \AND \bauthor{\bsnm{Vidakovic},~\bfnm{B.}\binits{B.}}
(\byear{1998}).
\btitle{Encyclopedia of Statistical Sciences}
\bvolume{6}.
\bpublisher{John Wiley \& Sons}.
\end{bbook}
\endbibitem

\bibitem[\protect\citeauthoryear{Stout}{1974}]{stout1974almost}
\begin{bbook}[author]
\bauthor{\bsnm{Stout},~\bfnm{William~F.}\binits{W.~F.}}
(\byear{1974}).
\btitle{Almost Sure Convergence}.
\bpublisher{Academic Press}.
\end{bbook}
\endbibitem

\bibitem[\protect\citeauthoryear{Zeileis et~al.}{2002}]{zeileis2002detecting}
\begin{barticle}[author]
\bauthor{\bsnm{Zeileis},~\bfnm{A.}\binits{A.}}, \bauthor{\bsnm{Leisch},~\bfnm{F.}\binits{F.}}, \bauthor{\bsnm{Hornik},~\bfnm{K.}\binits{K.}} \AND \bauthor{\bsnm{Kleiber},~\bfnm{C.}\binits{C.}}
(\byear{2002}).
\btitle{Detecting Structural Breaks in Multivariate Time Series}.
\bjournal{Journal of Time Series Analysis}
\bvolume{23}
\bpages{545--567}.
\end{barticle}
\endbibitem

\end{thebibliography}
\end{document}